\documentclass[journal]{IEEEtran}

\IEEEoverridecommandlockouts                              

\usepackage{color}
\usepackage{multirow}
\usepackage{cite}
\usepackage{balance}
\usepackage{enumitem}
\usepackage{color}
\usepackage{multirow}
\usepackage{cite}
\usepackage{balance}
\usepackage{amssymb}
\usepackage{soul}
\usepackage{pifont}
\usepackage{url}
\usepackage{tabularx}
\usepackage{array}
\newcolumntype{L}{>{\centering\arraybackslash}m{1.0cm}}
\usepackage[ruled,linesnumbered]{algorithm2e}

\usepackage[pdftex]{graphicx}
\usepackage[table,xcdraw]{xcolor}
\DeclareGraphicsExtensions{.pdf,.jpeg,.png,.JPG,}
\usepackage{graphicx}

\DeclareGraphicsExtensions{.eps}
\DeclareGraphicsExtensions{.tif}

\graphicspath{{figs/}}
\usepackage{epstopdf}
\usepackage{booktabs}
\usepackage{comment}

\usepackage{amsmath}
\usepackage{cleveref}
\ifCLASSOPTIONcompsoc
    \usepackage[caption=false, font=normalsize, labelfont=sf, textfont=sf]{subfig}
\else
\usepackage[caption=false, font=footnotesize]{subfig}
\fi

\begin{document}

\title{Network-Level Optimization for Unbalanced Power Distribution System: Approximation and Relaxation}

\author{Rahul Ranjan Jha,~\IEEEmembership{Student Member,~IEEE}, and~Anamika~Dubey,~\IEEEmembership{Member,~IEEE}
\thanks{R. Jha  and A. Dubey are with the School of Electrical Engineering and Computer Science, Washington State University, Pullman, WA, 99164, e-mail: rahul.jha@wsu.edu, anamika.dubey@wsu.edu.}}
\maketitle

\begin{abstract}
The nonlinear programming (NLP) problem to solve distribution-level optimal power flow (D-OPF) poses convergence issues and does not scale well for unbalanced distribution systems. The existing scalable D-OPF algorithms either use approximations that are not valid for an unbalanced power distribution system, or apply relaxation techniques to the nonlinear power flow equations that do not guarantee a feasible power flow solution. In this paper, we propose scalable D-OPF algorithms that simultaneously achieve optimal and feasible solutions by solving multiple iterations of approximate, or relaxed, D-OPF subproblems of low complexity. The first algorithm is based on a successive linear approximation of the nonlinear power flow equations around the current operating point, where the D-OPF solution is obtained by solving multiple iterations of a linear programming (LP) problem. The second algorithm is based on the relaxation of the nonlinear power flow equations as conic constraints together with directional constraints, which achieves optimal and feasible solutions over multiple iterations of a second-order cone programming (SOCP) problem. It is demonstrated that the proposed algorithms are able to reach an optimal and feasible solution while significantly reducing the computation time as compared to an equivalent NLP D-OPF model for the same distribution system.
  
\end{abstract}


\IEEEpeerreviewmaketitle

\section*{Nomenclature}
	\addcontentsline{toc}{section}{Nomenclature}
	\begin{IEEEdescription}[\IEEEusemathlabelsep\IEEEsetlabelwidth{$V_1,V_2$}]
        \item [Sets]
        \item[$\mathcal{G}$] Directed graph for distribution system
        \item [$\mathcal{E}$] Set of edges (branches) in $G$
        \item [$\mathcal{N}$] Set of buses (nodes) in $G$		
		\item[$\mathcal{N_{DG}}$] Set of nodes with smart inverter connected DGs
        \item[$\Phi_i$] Set of phases of bus $i$ where, $\Phi_i \subseteq \{a,b,c\}$
        \item[$\Phi_{ij}$] $\{(pq): p\in \Phi_i, q\in \Phi_j\}$
    \item [$f(x)$] Represents a general objective function
    \item [$g(x)$] Represents  a general linear and non-linear sets of equations
        \vspace{0.3cm}
        \item [Variables and Parameters]
        \item[$e^{pp}_{ij}$, $e^{pq}_{ij}$ ] Feasibility gap defined as error 
        \item[$I_{ij}^p$] $I_{ij}^p = |I_{ij}^p|\angle{\delta_{ij}^p}$ is complex line current corresponding to phase $p\in\Phi_{i}$ where, $|I_{ij}^p|$ is magnitude and $\delta_{ij}^p$ is corresponding phase angle 
        \item[$l_{ij}^p$] $l_{ij}^p= I_{ij}^p(I^p_{ij})^H$.
         \item[$I_{ij}^{pq}$] $I_{ij}^{pq} = |I_{ij}^{pp}| \times |I_{ij}^{qq}|$ product of  phase current where, $(pq)\in\Phi_{ij}$
         \item[$l$, $u$] Lower and upper bounds for the control variables
         \item[$m_i$, $n_i$] Non-negative constraints on penalty function
         \item[$p_{DG,i}^p$] Per-phase active power generated by DG connected at $i^{th}$ bus
		\item[$q_{DG,i}^p$] Available per-phase reactive power from DG connected at $i^{th}$ bus
        \item [$s$] Substation bus where, $s\in\mathcal{N}$
        \item[$s_{DG,i}^{p}$] $s_{DG,i}^p = p_{DG,i}^p + j q_{DG,i}^p$ per-phase is complex power for DG connected at bus $i \in \mathcal{N}$
        \item[$s_{DG,i}^{rated,p}$] Rated per-phase apparent power capacity for DG connected at bus $i \in \mathcal{N}$
        \item[$s_{L,i}^p$] $s_{L,i}^p = p_{L,i}^p + j q_{L,i}^p$ is complex power demand at bus $i$ corresponding to phase $p\in\Phi_i$ where, $p_{L,i}^p$ and $q_{L,i}^p$ are corresponding active and reactive power demand, respectively.
        \item[$S_{ij}^{pq}$] $S_{ij}^{pq} = P_{ij}^{pq} + j Q_{ij}^{pq}$ is complex power flow in branch $(i,j)$ corresponding to $(pq)\in\Phi_{ij}$, where, $P_{ij}^{pq}$ and $Q_{ij}^{pq}$ are corresponding active and reactive components, respectively.
        \item[$s^{k}$] Step bounds on the control variables
        \item[$V_i^p$, $v_{i}^p$] $V_{i}^p = |V_{i}^p|\angle{\theta_{i}^p}$ is complex voltage for $p\in\Phi_i$ and $v_{i}^p = |V_{i}^p|^2$ is square of voltage for $p\in\Phi_i$
        \item[$\mathcal{W}$] Weight of penalty function
        \item[$z_{ij}$] Complex three-phase impedance matrix for an edge $(i,j)\in \mathcal{E}$
		\item[$z_{ij}^{pq}$] $z_{ij}^{pq} = r_{ij}^{pq} + j x_{ij}^{pq}$ is an element of the complex impedance matrix $z_{ij}$ for a branch $(i,j)$ where, $(pq)\in\Phi_{ij}$
\end{IEEEdescription}

\section{Introduction}
\label{Introduction}

In an active power distribution system, optimal power flow (OPF) algorithms find multiple applications, including, but not limited to, loss minimization, volt-var optimization, and effective management of  distributed energy resources (DERs) \cite{molzahn2019survey}. The increasing penetration of DERs, the proliferation of proactive loads, and the interest in demand response programs require optimization methods for large-scale power distribution grid \cite{dubey2020paving}. Lately, these new requirements have encouraged the rapid adoption of advanced distribution management systems (ADMS) and related ADMS applications \cite{ngo2019investing}. As the distribution systems continue to become more active, the need for faster management of the grid's controllable assets will inevitably necessitate faster OPF algorithms \cite{sensitivitybased}.  Compared to the bulk power grid, distribution-level OPF (D-OPF) poses unique challenges due to three-phase unbalanced loading, mutual coupling among the different phases of the line, the existence of single-phase and two-phase branches, and radial topology with high R/X ratio leading to significant voltage drops.  While earlier work focused on the balanced distribution systems\cite{Low1,baran1989optimal}, lately, significant strides have been made regarding three-phase unbalanced D-OPF formulations. Unfortunately, these models lead to non-convex optimization problems, which are difficult to solve, especially for a large distribution system \cite{Spaudyal,jha2019bi, Bruno2011}.

Lately, several scalable D-OPF algorithms have been proposed in the literature, broadly categorized as methods based on (1) linear approximations; and (2) convex relaxations. A detailed survey on different approaches to solve  D-OPF  is presented in \cite{molzahn2019survey}. Linear approximations scale well for large systems \cite{ bernstein2017linear,  yang2018linearized, Bolognani}; but they either ignore power losses, or assume balanced system conditions, or assume node voltages to be close to their nominal values (1 pu). These assumptions are not valid for unbalanced power distribution systems that typically observe unbalanced loading conditions with significant power losses and voltage drops \cite{Dhople2015}. Further, the solutions of linearized D-OPF models are typically not feasible for the original nonlinear programming (NLP) OPF problem. One approach is to employ successive linear programming (SLP) methods where the basic idea is to solve the NLP OPF as multiple iterations of approximate linear programming (LP) problems. This simultaneously leads to a feasible and optimal solution. This approach has been explored to solve OPF for the bulk power grid \cite{ castillo2016successive,watson2014current}. However, the existence of mutual coupling among the phases and the requirement for solving OPF on the full three-phase model makes it challenging to apply the SLP algorithm for solving the D-OPF problems.  It is to be noted that while the exact linearization of three-phase power flow equations, as proposed in \cite{Bernstein2018}, may be used to develop SLP problems, the resulting linearized power models are computationally inefficient for large feeders. 
 To address this concern, authors in \cite{Bolognani2015, Bernstein2018} proposed a linear approximate power flow method using fixed point linearization (FPL).
However, the FPL method leads to significant bus voltage magnitude errors that are in the order of 1-1.5\%. 

Another approach attempts to relax the nonlinear power flow equations as convex inequalities. This results in a convex optimization problem for D-OPF, which solves within a reasonable time for large-scale distribution systems \cite{ gan2014convex,zhao2017convex}. These methods either relax a bus-injection model (BIM) based D-OPF as a semi-definite program (SDP) \cite{XBai}, or a branch-flow model (BFM) based D-OPF as a second-order cone program (SOCP)  \cite{ jabr2006radial}. Although the proposed SDP and SOCP relaxations lead to a convex problem of reduced complexity, they may result in solutions that are infeasible for the original nonlinear power flow model. Consequently, several researchers have attempted to derive conditions that ensure the exactness of the relaxed OPF problems \cite{ jabr2006radial, MFarivar, Low2}. While convex relaxations for single-phase distribution systems have been found to be exact for a certain choice of objective functions and under specific conditions on the distribution system, no such guarantees exist for a three-phase unbalanced system \cite{Nanpeng}. For example, SDP relaxation has been applied to three-phase D-OPF problems \cite{gan2014convex, DallAnese2013Distributed}. However, it has been reported in recent articles that for three-phase distribution systems the SDP relaxed model may lead to numerical stability issues \cite{gan2014convex, wang2017accurate}, and infeasible power flow solutions depending upon the choice of system parameters and objective function \cite{DallAnese2018, Nanpeng,jhaitersocp}. 
Thus, in all these scenarios, new algorithms are needed that can still make use of convex optimization techniques that are computationally attractive but can simultaneously result in feasible power flow solutions. In \cite{Nanpeng}, the authors proposed a convex iteration technique to solve SDP relaxation for the D-OPF model that leads to a feasible and optimal solution. However, the approach does not scale well for large systems and requires additional heuristics. In \cite{DallAnese2018}, authors proposed an iterative approach that starts with a feasible operating point for power flow and solves multiple iterates of convex programming problems to reach to the optimal solution. This algorithm takes 1000 iterations for the IEEE 37-bus multi-phase test system, where each iteration takes 2-5 sec (total 2000-5000 sec $\approx$ 33-83 min) with a flat start. However, with warm start (starting with close enough voltage profile from the last time-step), it converges within 20 sec. Notice that for high-levels of DER penetrations, the voltage profiles can significantly vary from one time-step to the next. Thus, the implications of warm start in speeding up the algorithm in a general setting and for larger feeder models needs further analysis. Similarly, in our previous work, we proposed an iterative approach where a feasible OPF solution is obtained by solving multiple iterations of relaxed-OPF problems; however, the approach is only applicable to a single-phase system \cite{jhaitersocp}. 

The major limitations posed by the existing literature on D-OPF methods are: (1) nonlinear OPF models do not scale well for large-scale distribution systems especially when considering phase unbalances and phase-frame network model; (2) linearized D-OPF models are based on approximations that are not valid for unbalanced distribution systems and lead to sub-optimal solutions; (3) SLP methods are not directly applicable for solving D-OPF problems given the mutual coupling among the phases that introduces additional nonlinearities; (4) convex relaxations techniques applied to D-OPF problems are not guaranteed to be exact with respect to original power flow model resulting in infeasible operating points; and  (5) the existing iterative methods that employ convex iteration techniques are computationally challenging for large multi-phase distribution feeders \cite{Nanpeng,DallAnese2018}. While the applications of distributed optimization techniques has also been extensively explored for solving D-OPF problems, they also pose the aforementioned limitations. Furthermore, they require hundreds, if not thousands, of communication rounds with the neighboring agents to reach a converged solution and pose convergence challenges, both practically and theoretically, for general nonlinear D-OPF problems \cite{ molzahn2017survey}.

In this paper, we propose two iterative algorithms based on linear approximation and convex relaxation to solve D-OPF problems within a reasonable time. The proposed scalable D-OPF formulations simultaneously achieve optimal and feasible solutions for unbalanced power distribution systems.
The results are demonstrated using IEEE 123-bus (with 267 single-phase nodes) and modified R3-12.47-2 feeder (with 860 single-phase nodes). The proposed algorithms are shown to converge to the solutions obtained upon solving the NLP D-OPF model while significantly improving the computation time. While the linear approximation approach converges faster, the method based on convex relaxation is useful in providing bounds on the true optima for the NLP problem. 
The contributions of the paper are listed below:
\begin{itemize}[leftmargin=*]
\item {\em Iterative D-OPF algorithms for unbalanced power distribution systems} – We propose iterative algorithms for three-phase AC OPF for radial power distribution systems that achieve optimal and feasible D-OPF solutions. The proposed approach is applicable when convex relaxation, such as Semi-definite relaxation, may fail to converge to a feasible power flow solution. The proposed methods use computationally efficient optimization models, LP and SOCP, to solve the difficult NLP problem for large feeders. 
\item {\em Scalable D-OPF algorithm using approximate linear power flow equations: } We develop a scalable D-OPF algorithm by formulating NLP D-OPF as a penalty-successive linear programming (PSLP) problem. The proposed approach formulates LP sub-problems using the first-order Taylor series approximation of the linear-quadratic power flow model introduced in \cite{jha2019bi} around the current iterate. The PSLP algorithm then solves multiple iteration of the LP problem via successive linearization. The algorithm converges within 3-4 iterations of the LP problem.  
\item {\em Scalable D-OPF algorithm using relaxed nonlinear power flow equations:} We propose a scalable D-OPF algorithm by formulating NLP D-OPF as an iterative second-order cone programming (ISOCP) problem. The proposed ISOCP algorithm solves multiple iterates of the SOCP problem to achieve optimal and feasible power flow solutions. Additional directional constraints are added to ensure that the feasibility gap is reduced over successive SOCP iterations. We also detail the convergence and optimality of the ISOCP algorithm. 
\end{itemize}

\section{Background: Three-phase Power Flow Model and Approximations}
We introduce the mathematical formulation for the three-phase branch flow model (BFM). The proposed approximations to derive a BFM of lower complexity are discussed next, followed by the resulting approximate BFM formulation.

\vspace{-0.2cm}
\subsection{Three-Phase Power Flow using Branch Flow Model}
A radial distribution system can be represented as a directed graph $\mathcal{G} = (\mathcal{N}, \mathcal{E})$ where $\mathcal{N}$ denotes set of buses and $\mathcal{E}$ denotes set of edges. An edge $(i,j)$ joins two adjacent nodes $i$ and $j$ where $i$ is the parent node for node $j$. The three phase $\{a,b,c\}$ for a node $i$ in the distribution system is denoted by  $\Phi_i \in \{a,b,c\} $.  For each bus $i \in \mathcal{N}$ and phase $p$, let $V_i^{p}$, $s_{L,i}^p$, and $s_{DG,i}^p$ be complex voltage, complex power demand and complex DG power generation, respectively. Let, $V_i := [V_i^{p}]_{p \in \Phi_i}$, $s_{L,i} := [s_{L,i}^{p}]_{p \in \Phi_i}$ and $s_{DG,i} := [s_{DG,i}^{p}]_{p \in \Phi_i}$. For each line, let $p$ phase current be $I_{ij}^{p}$ and define, $I_{ij} := [I_{ij}^{p}]_{p \in \Phi_i}$. Let $z_{ij}$ be the phase impedance matrix.

The mathematical formulation for a power flow model based on branch flow equations for a radial distribution system is detailed in (\ref{eq1})-(\ref{eq4}) \cite{gan2014convex}. For sake of simplicity, here we have removed the phase notation.   The voltage drop and power balance equations are given by (\ref{eq1}) and (2), respectively. The relationship between the branch power flow, nodal voltages, and branch currents is defined using (3) and (4). Note that $(.)^H$ represents the conjugate transpose. 

Note that a D-OPF model based on the power flow equations in (1)-(4) is non-convex due to the  rank-1 constraint for the positive semidefinite (PSD) matrix in (4). While the existing body of literature presents several methods to obtain a convex D-OPF formulation, the relaxed formulations are not guaranteed to be exact for an unbalanced three-phase power distribution system \cite{Nanpeng,DallAnese2018}. Furthermore, relaxing the rank-1 constraint leads to a SDP model that does not scale well. This calls for new formulations for power flow model and new algorithms for solving OPF for an unbalanced power distribution system. 

\vspace{-0.3cm}
\begin{small}
\begin{eqnarray}\label{eq1}
v_i  = v_j + (S_{ij}z_{ij}^H+z_{ij}S_{ij}^H) - z_{ij}l_{ij} z_{ij}^H \\
\text{diag}(S_{ij} - z_{ij}l_{ij}) - s_{L,j} + s_{DG,i}  =  \sum_{k:j \rightarrow k}{\text{diag}(S_{jk})}   \\
\left[
  \begin{array}{cc}
    v_i & S_{ij} \\
    S_{ij}^H & l_{ij} \\
  \end{array}
\right] = \left[
  \begin{array}{c}
    V_i\\
    I_{ij}\\
  \end{array}
\right]
\left[
  \begin{array}{cc}
    V_i\\
    I_{ij}\\
  \end{array}
\right]^H \\
\left[
  \begin{array}{cc}\label{eq4}
    v_i & S_{ij} \\
    S_{ij}^H & l_{ij} \\
  \end{array}
\right] : - \text{Rank-1 PSD Matrix}
\end{eqnarray}
\end{small}
\vspace{-0.2cm}

\vspace{-0.3cm}
\subsection{Approximate Three-Phase Power Flow Equations}

Our objective is to reduce the complexity of the branch flow equations in (1)-(4). Towards this goal, we introduce valid approximations for power flow variables that help us avoid the matrix-based formulation in (\ref{eq1})-(\ref{eq4}). This results in a new three-phase BFM where power flow equations are expressed as a set of linear and quadratic equations \cite{jha2019bi}. The novelty lies in carefully constructed approximations on voltage and current angles that help simplify the phase-coupled voltage and power balance equations into a per-phase coupled description. 

Define, for
\begin{small}
$p \in  \Phi_{i}$: $V_{i}^p = |V_{i}^p|\angle{\theta_{i}^p}$, $s_{L,i}^p = p_{L,i}^p + j q_{L,i}^p$, $s_{DG,i}^p = p_{DG,i}^p + j q_{DG,i}^p$, $v_i^p = (V_i^p)^2$; for $(pq) \in  \Phi_{ij}$: $l_{ij}^{pp} = (I_{ij}^{pp})^2$, $I_{ij}^p = |I_{ij}^p|\angle{\delta_{ij}^p}$, $l_{ij}^{pq} = (|I_{ij}^p|\times |I_{ij}^q|)$, $\delta_{ij}^{pq} = \delta_{ij}^{p} - \delta_{ij}^{q}$, $S_{ij}^{pq} = P_{ij}^{pq} + j Q_{ij}^{pq}$, and $z_{i,j}^{pq} = r_{ij}^{pq} + j x_{ij}^{pq}$, 
$ (l_{ij}^{pq})^2 = l_{ij}^{pp} \times l_{ij}^{qq}$.  
\end{small}

\vspace{0.1cm}
\begin{itemize}[nolistsep,leftmargin=*]
\item {\em Approximating Nodal Voltage Phase Angle:} For a given node, it is assumed that the nodal voltage phase angles are separated by $120^0$ and the degree of unbalance in voltage magnitudes is not large. This assumption allows us to represent off diagonal elements $S_{ij}^{pq}$ as a function of the diagonal elements, $S_{ij}^{pp}$, in $S_{ij}$. Note that the allowable limits on bus voltage magnitudes and phase unbalance as per the ANSI C83.1 standard justifies this assumption. Additional validation of this assumption is presented using different test feeder models in our prior work \cite{jha2019bi}.

\item {\em Approximating Angle Difference between Phase Currents:} On expanding (1) and (2), nonlinearities are introduced as trigonometric functions of angle difference between the phase currents that significantly increase the complexity of the OPF problem. 
That is, we observe terms corresponding to $\sin(\delta_{ij}^{pq})$ and $\cos(\delta_{ij}^{pq})$ in power flow expressions, where, $\delta_{ij}^{pq} = \angle{\delta_{ij}^p} - \angle{\delta_{ij}^q}$. In the proposed formulation, the phase angle differences between branch currents are approximated and modeled as constant variables and equal to the one obtained by solving power flow with constant impedance load model. Note that the constant impedance load model is only used to approximate $\delta_{ij}^{pq}$ and this assumption does not limit the type of load that can be incorporated in the proposed OPF model, see \cite{jha2019bi} for details. 
\end{itemize}

Using these assumptions and upon expanding (\ref{eq1})-(\ref{eq4}) on per-phase basis, we obtained a set of linear and quadratic equations detailed in (\ref{eq6})-(\ref{eq10}). Here, (\ref{eq6}) and (6)  are active and reactive power balance equations for the phase $pp \in \Phi_{ij}$ of branch $(i,j)$, respectively; (\ref{eq8}) is the voltage drop equation for  phase $p \in \Phi_i$ of nodes $i$ and $j$; (\ref{eq9}) represents the relationship between the nodal voltage, branch current, and branch power flow for the phase $pp \in \Phi_{ij}$ of branch $(i,j)$ and phase $p \in \Phi_i$ of  nodes $i$ and $j$; (\ref{eq10}) represents the relation between the branch current in phases $pp \in \Phi_{ij}$ and $qq \in \Phi_{ij}$ of branch $(i,j)$. Here, (\ref{eq6})-(\ref{eq8}) are linear while (\ref{eq9})-(\ref{eq10}) are quadratic wrt. D-OPF variables, 
 \begin{small}
 $x_{OPF}$: $\{q_{DG,i}^p$, $v_i^p$, $P_{ij}^{pp}$, $Q_{ij}^{pp}$, $P_{ij}^{qq}$, $Q_{ij}^{qq}$, $S_{ij}^{pq}$, $l_{ij}^{pp}$, $l_{ij}^{qq}$, $l_{ij}^{pq}$, $\forall (i,j)\in\mathcal{E}, \forall i\in\mathcal{N}, \forall p \in \Phi_i, \forall pp \in \Phi_{ij}$\}. 
\end{small}

 \vspace{-0.1cm}
 \begin{small}
 \begin{eqnarray}\label{eq6}
  \nonumber &P_{ij}^{pp} - \sum_{q \in \Phi_j}{l_{ij}^{pq}\left(r_{ij}^{pq} \cos(\delta_{ij}^{pq})- x_{ij}^{pq} \sin(\delta_{ij}^{pq})\right)}  \\ &= \sum_{k:j \rightarrow k}P_{jk}^{pp} + p_{L,j}^p - p_{DG,j}^p \\
  \nonumber&Q_{ij}^{pp} - \sum_{q \in \Phi_j}{l_{ij}^{pq}\left(x_{ij}^{pq} \cos(\delta_{ij}^{pq})+ r_{ij}^{pq} \sin(\delta_{ij}^{pq})\right)} \\&= \sum_{k:j \rightarrow k}Q_{jk}^{pp} + q_{L,j}^p - q_{DG,j}^p 
\end{eqnarray}
\vspace{-0.5cm}
\begin{eqnarray}\label{eq8}
 \nonumber &v_i^p = v_j^p + \sum_{q \in \Phi_j}{2 \mathbb{Re}\left[S_{ij}^{pq} (z_{ij}^{pq})^*\right]} - \sum_{q \in \Phi_j}{(z_{ij}^{pq})^2 l_{ij}^{qq}} \\
  &- \sum_{q1,q2 \in \Phi_j, q1 \neq q2}{2\mathbb{Re}\left[ z_{ij}^{pq1} l_{ij}^{q1q2}\left(\angle(\delta_{ij}^{q1q2})\right)(z_{ij}^{pq2})^*\right]}
\end{eqnarray}
\vspace{-0.5cm}
\begin{eqnarray}\label{eq9}
  (P_{ij}^{pp})^2 + (Q_{ij}^{pp})^2 = v_i^p  l_{ij}^{pp}
\end{eqnarray}
\vspace{-0.5cm}
\begin{eqnarray}\label{eq10}
  (l_{ij}^{pq})^2 = l_{ij}^{pp}  l_{ij}^{qq}
\end{eqnarray}
\end{small}
\vspace{-0.5cm}
 
 Thus, the first approximation of phase voltage angles for each allows us to convert the matrix-based formulation in (1)-(4) to the sets of equations in (5)-(9). Then, replacing the current angles with pre-calculated values renders the trigonometric terms in (5)-(9) as constants and results in a linear-quadratic power flow model.

We have thoroughly validated the above two assumptions in our prior work for several distribution feeder models with varying degree of voltage unbalances and for different load models \cite{jha2019bi}. Note that in the proposed approximate power flow model, the number of variables has been reduced from $36 \times(n-1)$ to $15\times(n-1)$, where $n$ is the number of nodes. While it is of lower complexity compared to (1)-(4), the formulation is still difficult to scale due to nonlinear quadratic equations (8) and (9). In what follows, we present new D-OPF algorithms based on the approximation and relaxation of the quadratic equality constraints in the approximate branch flow model.

\noindent \textbf{Discussion - }{\em Justification for Linear-quadratic model in the context of state-of-art methods:} The first-order approximation of the exact three-phase power flow model, (1)-(3), as developed in \cite{Bernstein2018} is relatively accurate; however, they are computationally expensive for large networks. On the other hand, other approximate models such as fixed-point linearization (FPL) \cite{Bolognani2015} and linearization around 1 pu voltage \cite{Dhople2015} leads to significant errors in approximating the node voltage magnitude. Since the OPF problems are generally challenging, and the integration of new devices will only add to the complexity, approximations to power flow equations will be inevitable. In this context, the proposed linear-quadratic model leads to a simpler first-order linearization that is significantly more accurate \cite{jha2019bi}. 

The aforementioned linear-quadratic approximation of power flow is motivated by the domain-specific knowledge related to the power flow characteristics in radial power distribution systems. We argue that the domain-specific knowledge can better inform the algorithms and the associated approximations. The use of domain-specific knowledge in developing approximate models in power systems is not new. For example, the DC power flow model that approximates node voltage magnitude to 1 pu works very well for the transmission systems; it fails for the distribution networks. On the contrary, in radial distribution networks, dist-flow equations provide a good approximation of AC power flow \cite{baran1989optimal}. Likewise, while the bulk power grid uses the Newton-Raphson method for power flow, the distribution power flow solvers commonly use a forward-backward sweep algorithm (fixed-point iteration) that leverages the radial distribution feeder topology. Thus, given the long history of domain-specific approximations in the power systems community, we believe such formulations can lead to more useful results for D-OPF problems.

\section{Optimal Power Flow for Three-Phase System}
In this section, we formulate the D-OPF problem for an unbalanced distribution system based on the approximate power flow model detailed in (5)-(9). 
The problem objective is to reduce the active power consumption for the distribution feeder by controlling the bus voltages.  For a distribution system with voltage dependent loads, the required objective can be achieved by operating all the nodes voltage towards their lower limits.  While network losses may increase upon decreasing the node voltages, the reduction is load demand (due to voltage-dependent loads) is typically larger than the increase in network losses, as observed in smart grid pilots \cite{AvistaCVR}. By properly coordinating the reactive power dispatch of the smart inverters of the DGs interfaced with the feeder, a reduction in node voltages can be achieved leading to a reduction in feeder active power consumption. This technique for  power  reduction is known as conservation voltage reduction (CVR).

The D-OPF problem to minimize the feeder active power consumption by coordinating the reactive power dispatch from smart inverters is modeled as an NLP problem detailed in (5)-(16). Notice that in the proposed model, the choice of objective function is arbitrary and can be modified to other problem objectives such as loss minimization, feeder voltage regulation.

\vspace{0.2cm}
\begin{minipage}{23.5em}
		\flushleft
		\small
		Variables: $x_{OPF}$: $\{q_{DG,j}^p$, $v_i^p$, $P_{ij}^{pp}$, $Q_{ij}^{pp}$, $P_{ij}^{qq}$, $Q_{ij}^{qq}$, $S_{ij}^{pq}$, $l_{ij}^{pp}$, $l_{ij}^{qq}$, $l_{ij}^{pq}$, $\forall (i,j)\in\mathcal{E}, \forall i\in\mathcal{N}, \forall p \in \Phi_i, \forall pp \in \Phi_{ij}$\}
		\begin{flalign}\label{eq15nl}
		\text{Minimize:}  \ \ \sum_{p\in \Phi_s,j:s \rightarrow j}{P_{sj}^p}		&&
		\end{flalign}
		Subject to:
        \begin{flalign}\label{eq17nl}
        \nonumber (5)-(9) \hspace{12pt} &\forall (i,j)\in\mathcal{E}, \forall i\in\mathcal{N}\\
        p_{L,i}^p = p_{i,0}^p + CVR_{p}\dfrac{p_{i,0}^p}{2}(v_i^p-1)  \hspace{12pt} &\forall i\in\mathcal{N_L} \\
        q_{L,i}^p = q_{i,0}^p + CVR_{q}\dfrac{q_{i,0}^p}{2}(v_i^p-1)  \hspace{12pt} &\forall i\in\mathcal{N_L} \\
        q_{DG,i}^p \leq \sqrt{(s_{DG,i}^{rated,p})^2 - (p_{DG,i}^p)^2} \hspace{12pt} &\forall (i) \in \mathcal{N_{DG}} \\
        q_{DG,i}^p \geq -\sqrt{(s_{DG,i}^{rated,p})^2 - (p_{DG,i}^p)^2} \hspace{12pt} &\forall (i) \in \mathcal{N_{DG}}\\
        (V_{min})^2\leq v_i^p \leq (V_{max})^2 \hspace{12pt} &\forall i\in \mathcal{N}\\
         l_{ij}^{pp} \leq (I^{rated}_{ij})^2 \hspace{12pt} &\forall i\in \mathcal{N}
        \end{flalign}
        \vspace{-0.2cm}
\end{minipage}

where, (10) states the problem objective of minimizing the active power consumption from the substation; constraints (5)-(9) are approximate nonlinear AC power flow equations; constraints (11)-(12) define CVR-based voltage dependent load model; constraints (13)-(14) define the limits on the reactive power dispatch from smart inverters, where, $s_{DG,i}^{rated,p}$ is the rated per-phase apparent power capacity for the DG at bus $i \in \mathcal{N_{DG}}$; constraint (15) limits the node voltages where, $V_{min} = 0.95$ and $V_{max} = 1.05$; and constraint (16) limits the branch current flow by the thermal rating of the lines $I^{rated}_{ij}$.

{\it Discussions: } A typical distribution feeder also includes legacy voltage control devices such as, capacitor banks and voltage regulators with discrete control variables. Thus, D-OPF problem is typically modeled as a mixed integer nonlinear programming (MINLP) problem. In our prior work, we proposed a two-stage formulation to solve the MINLP problem by decomposing the problem into an MILP (Stage-1) and NLP (Stage-2) \cite{jha2019bi}. While Stage-1 (MILP) scales well for larger feeders, Stage-2 (NLP) did not. Thus, in this paper, our objective is to improve the scalability of the associated nonlinear programming (NLP) problem for D-OPF formulations. We, therefore, ignore the integer variables for the discussion herein. 

\section{Proposed Iterative OPF Algorithms using Approximation and Relaxation Techniques}
The D-OPF described in (5)-(16), while can be  solved using non-linear optimization solvers, the computation time significantly increases with the increases in the problem size. Further, it gets more difficult to achieve convergence for larger problem size. The complexity arise due to the nonlinear equality constraints in the BFM, i.e. (8)-(9). In this paper, we propose two approaches to reduce the problem complexity and achieve superior convergence and computational speed. These methods are based on the approximations and relaxations of the nonlinear equality constraints, (8)-(9), and innovations to drive the approximate and relaxed solutions towards the feasible space for the NLP D-OPF problem. 
\begin{itemize}[leftmargin=*]
 \item Approximation - In this approach, the nonlinear equality constraints are linearized around the current operating point. This results in a LP problem that is iteratively solved to achieve an optimal and feasible solution. However, it is possible that when approximating we exclude portions of the original feasible space and the formulation may converge to sub-optimal points.   
 \item Relaxation - In this approach, the nonlinear equality constraints are relaxed as inequalities resulting in conic constraints. This reduces the NLP D-OPF to a second-order cone programming (SOCP) problem. The SOCP problem is a convex problem that is easier to solve and scale well for large problem size. The relaxation, however, leads to an extension of the feasible space. Thus, the resulting solution may not be feasible wrt. to the original NLP. Here, we add additional directional constraint to drive the solution to feasible space over multiple iterations of SOCP problem. 
\end{itemize}

\subsection{Approximation -- D-OPF via Iterative Linear Programming}
The proposed approach is based on the successive linearization of the nonlinear equality constraints in BFM described in (5)-(9) around the current operating point, and iteratively solving the resulting LP problem until the optimal solution is obtained. In what follows, first, we describe the approach to linearize the nonlinear equality constraints in (8)-(9). Next, we discuss the proposed penalty successive linear programming (PSLP) algorithm to solve the original NLP D-OPF problem. 

\subsubsection{Linearized Power Flow Model - Taylor Series Approximation}
The power flow equations (5)-(7) are the linear equations wrt. problem variables, denoted by the set $x_{OPF}$. However, (8)-(9) are nonlinear equations wrt. the power flow variables. It is known that any nonlinear equations can be linearized around an operating point using the first-order Taylor series approximation. Here, we linearize (8)-(9) around the current operating point denoted by $P_{ij,0}^{pp}$, $Q_{ij,0}^{pp}$, $v_{i,0}^p$,  and  $v_{i,0}^q$. The linearized equations are detailed in (\ref{eqlpp}) and (\ref{eqlpq}).

\begin{equation}\label{eqlpp}
\small
  \left(l_{ij}^{pp}\right) =   \dfrac{2(P_{ij,0}^{pp}P_{ij}^{pp}+Q_{ij,0}^{pp}Q_{ij}^{pp}) }{(v_{i,0}^{p})}- \dfrac{(S_{ij,0}^{pp})^2 v_{i}^{p}}{(v_{i,0}^{p})^2}
 \end{equation} 

\begin{small}
 \begin{flalign} \label{eqlpq}
 \left(l_{ij}^{pq}\right) & =  \dfrac{1}{2}  \left( \dfrac{S_{ij,0}^{qq}}{S_{ij,0}^{pp}}\sqrt{\left(\dfrac{v_{i,0}^p}{v_{i,0}^q}\right)} l_{ij}^{pp} + \dfrac{S_{ij,0}^{pp}}{S_{ij,0}^{qq}}\sqrt{\left(\dfrac{v_{i,0}^q}{v_{i,0}^p}\right)}l_{ij}^{qq}\right) 
 \end{flalign}  
 \end{small}

\noindent where, {\small $S_{ij,0}^{pp} = \sqrt{(P_{ij,0}^{pp})^2 + (Q_{ij,0}^{pp})^2}$}. The $l_{ij}^{pp}$ and  $l_{ij}^{qq}$ in (\ref{eqlpq}) can be obtained using (\ref{eqlpp}).  Notice that $P_{ij,0}^{pp}$, $Q_{ij,0}^{pp}$, $v_{i,0}^p$,  and  $v_{i,0}^q$ are constants that are the known operating points around which (8) and (9) are linearized. Thus, (\ref{eqlpp}) and  (\ref{eqlpq}) are linear equations wrt. to the problem variables, i.e. $x_{OPF}$. Following the linearization, we develop an iterative approach to successively solve the D-OPF problem based on the linearized power flow equations. Specifically, we employ successive linear programming (SLP) approach to solve the original NLP model for D-OPF.

\subsubsection{Penalty Successive Linear Programming (PSLP)} 
We formulate the D-OPF problem as a penalty successive linear programming (PSLP). In PSLP, a sequence of linear programming (LP) problems are solved by approximating the original nonlinear optimization problem around the current operating point. The iterations continue until the change in objective function value over two successive iterations are within a pre-defined tolerance. It is known that the linearization of a nonlinear model is accurate only in the close vicinity of the operating points that were used for the linearization. That is, a large step change over successive LP iterations may lead to inaccuracies as the linearization may not be valid. This may result in oscillations in decision variables over successive LP iterations. To avoid this, in PSLP, we impose bounds on the step change in decision variables over successive LP iterations.  Briefly, the convergence result states that the PSLP algorithm converges to a stationary point of the L1 exact penalty function of the original NLP problem. The proof is based on two theorems. First theorem states that {\em ``the changes predicted by the PSLP algorithm for the piecewise linear approximation of the penalty problem (for the original NLP) are non-negative, implying that the ratio of actual to predicted change is well defined".} The second theorem states that {\em ``if the constrained level set of the exact penalty function (at the initial point of the algorithm) is bounded and the sequence of iterates generated by the PSLP algorithm is infinite, then the sequence of iterates has limit points, and every limit point is a constrained stationary point of the exact penalty function of the original NLP".} Authors then also state that in most cases, such a stationary point is also a KKT point of the NLP.\\
\vspace{-0.1cm}
\begin{algorithm}[h]
\small
\caption{\small Three-phase OPF using Penalty Successive Linear Programming (PSLP)}\label{alg:PSLP}
\SetAlgoLined
    \SetKwInOut{Input}{input}
    \SetKwInOut{Output}{output}
\SetKwInOut{defi}{define}
\SetKwInOut{init}{initialize}
\defi{$x = x_{OPF}$}
\Input{distribution system connectivity model, impedance matrix ($z_{ij}$), load parameters ($p_{L,i}^p$,$q_{L,i}^p$), DG parameters ($p_{DG,i}^p$, $s_{DG,i}^{rated,p}$), operating limit ($I_{ij}^{rated}$), $tol = 0.0001$}
\Output{Reactive power dispatch from DGs ($q_{DG,i}^p$)}
\init{iteration count, $k = 1$; $\varepsilon^{(k)} = 1$; $s^{(k)} = 0.01$; $x^{(k)} = x_{lin}$, where, $x_{lin}$ is the solution of three-phase OPF using LinDistFlow.} 
\While {$|\varepsilon^{(k)}| > tol$}{
{Linearize (8)-(9) around $x^{(k)}$ using (17)-(18)}\\
{Solve PSLP equations in (19)-(23) with linearized power flow equations for $\Delta x^{(k+1)}$} \\
{Update problem variables: $x^{k+1} = x^k + \Delta x^{(k+1)}$}\\
{Calculate: $\varepsilon^{(k+1)} = f(x^{(k)}) - f(x^{(k+1)})$}\\
	\If {$\varepsilon^{(k+1)}>0$}
	{
     $s^{(k+1)} = \dfrac{s^{(k)}}{2}$
         \\
        \Else {$s^{(k+1)} = 2{s^{(k)}}$}
	}	
{Increment iteration count, $ k = k+1$} \\
}
\label{algo}
\vspace{-0.1cm}
\end{algorithm}
\vspace{-0.2cm}

The notional problem description for PSLP-based OPF is described next in (19)-(23). Let $f(x)$ represent the objective function and $g(x) = 0$ represents the set of linear and nonlinear power flow equations. Let $x^k$ be the value of $x$ at the $k^{th}$ iteration and $\Delta x^{(k+1)}$ is the change of variable. The $(k+1)^{th}$ iteration of the PSLP problem solves for the change in decision variable $\Delta x^{(k+1)}$ by solving the following linear programming (LP) problem in (19)-(23). 

\begin{minipage}{23.5em}
		\flushleft
		\small
		 \begin{flalign}\label{eqpslpl}
		\text{Minimize:}  \ \ f(x^k) + \nabla f \Delta x^{(k+1)} + 	\mathcal{W} \sum_{i}  (m_i + n_i)	&&
		\end{flalign}
		\vspace{-0.1 cm}
		Subject to:
		 \vspace{-0.5 cm}
        \begin{flalign}\label{eqpslplc}
        g_i(x^k) + \nabla g_i \Delta x^{(k+1)} -b_i = m_i - n_i \\
        -s^{(k+1)} \leq \Delta x^{(k+1)} \leq s^{(k+1)}  \\
        l \leq x^k + \Delta x^{(k+1)} \leq u  \\
        m_i \geq 0  \hspace{0.1 cm}\text{and} \hspace{0.1 cm} n_i \geq 0
         \end{flalign}
         \vspace{-0.3cm}
\end{minipage}

where, (19) is the augmented objective with penalty weight $\mathcal{W}$; (20) is the linearized AC power flow equations defined at $k^{th}$ iteration; (21) imposes step bound on the variables at $(k+1)^{th}$ iteration; (22) ensures that the new operating point satisfy the original bound; (23) is the non-negativity constraint on penalty coefficients, $m_i$, $n_i$. Notice the objective function and constraint sets are linearized around $x^k$ where, $x^k$ is the solution from the $k^{th}$ iteration of the LP problem. Upon solving for $(k+1)^{th}$ iteration, a new operating point is obtained as $x^{k+1} = x^k + \Delta x^{(k+1)}$.

The algorithm for PSLP-based OPF problem is detailed in Algorithm 1.
First, we obtain the initial operating point around which the iterations of the PSLP problem are defined. In this paper, the initial operating points are obtained by solving the three-phase LinDistFlow equations \cite{gan2014convex}. The nonlinear power flow equations are then linearized around the current operating point to obtain  (\ref{eqlpp})-(\ref{eqlpq}). At $(k+1)^{th}$ iteration, the PSLP is solved for $\Delta x^{(k+1)}$ and a new operating point is obtained using $x^{k+1} = x^k + \Delta x^{(k+1)}$. The difference in the objective function values at $k$ and $(k+1)^{th}$ iterations, i.e., $\varepsilon^{(k+1)} = f(x^{k}) - f(x^{k+1})$, is calculated and termed as the error variable, $\varepsilon^{(k+1)}$. If $\varepsilon^{(k+1)}$ is less than the set tolerance value, the algorithm stops; else, iteration continues. The bounds on the variable change, $s^{(k+1)}$, is adjusted during iterations to improve convergence.  Briefly, if the error, $\varepsilon^{(k+1)}>0$, we need to decrease the trust region, $s^{(k+1)}$, to reach towards optimal solution and vice-versa. The iterations continue until $\varepsilon^{(k+1)}$ is less than the pre-specified tolerance.
 
 \vspace{-0.4cm}
 \subsection{Relaxation -- OPF using Iterative Convex Programming}
In this section, we derive a convex programming model for the three-phase OPF problem by relaxing the quadratic equality constraints, (8)-(9) as conic constraints shown in (\ref{eqcone1})-(\ref{eqcone2}). This reduces the original NLP model for D-OPF to a SOCP problem. Note that relaxation increases the feasible space of the actual NLP D-OPF problem.
\begin{equation}\label{eqcone1}
\small
  (P_{ij}^{pp})^2 + (Q_{ij}^{pp})^2 <= v_i^{p} l_{ij}^{pp}
\end{equation}
\vspace{-0.4cm}
\begin{equation}\label{eqcone2}
\small
  (l_{ij}^{pq})^2 <= l_{ij}^{pp}  l_{ij}^{qq}
\end{equation}
\vspace{-0.3cm}

As with any relaxation-based formulations, it needs to be determined whether the solutions obtained using the relaxed model are feasible wrt. the original problem. Note that the D-OPF solutions obtained using SOCP model will be exact w.r.t. the NLP model if, and only if, the SOCP solution satisfies the quadratic equality constraints in the NLP D-OPF model, i.e., (8)-(9). While SOCP relaxations have been thoroughly validated for single-phase radial distribution systems and have been found to be exact for minimization problems, no such guarantee exists for three-phase OPF problems \cite{Low1}. In fact, for three-phase unbalanced systems, we have verified that the SOCP relaxation proposed using (24)-(25) is not exact for the original NLP problem. That is, the D-OPF solution obtained using SOCP is not feasible wrt. the NLP D-OPF problem. 

In this paper, we propose an iterative SOCP algorithm to obtain a feasible and optimal D-OPF solution by iteratively solving relaxed OPF problem i.e. the SOCP model. The proposed approach derives directional constraints to reduce the feasibility gap over successive SOCP iterations. In the following subsections, first, we derive the equations for directional constraints used to reduce the feasibility gap. Next, we detail the algorithm for the iterative SOCP algorithm to solve the D-OPF problem for an optimal and feasible solution. 

\subsubsection{Deriving Directional Constraints}
The D-OPF solutions using the iterative approach will be feasible w.r.t. the original D-OPF problem (before relaxation) only if: (1) the difference between the $(P_{ij}^{pp})^2 + (Q_{ij}^{pp})^2 $ and $v_i^pl_{ij}^{pp}$ in (24) is gradually reduced to zero; and (2) the difference between $(l_{ij}^{pq})^2$ and $l_{ij}^{pp}l_{ij}^{qq}$ in (25) is gradually reduced to zero. The proposed method is designed to specifically achieve this equality over successive SOCP iterations. 

We define error terms, $e^{pp(k)}_{ij}$ and $e^{pq(k)}_{ij}$, measuring the feasibility gap at $k^{th}$ iteration, defined as (26) and (27).

\vspace{-0.2cm}
\begin{small}
\begin{eqnarray}\label{fgap1}
 e^{pp(k)}_{ij}  &=& (P_{ij}^{{pp}(k)})^2 + (Q_{ij}^{pp(k)})^2 - v_i^{p(k)} l_{ij}^{pp(k)} \hspace{0.1cm}\forall {ij} \in \mathcal{E}  \\
 e^{pq(k)}_{ij} &=& (l_{ij}^{pq(k)})^2  - l_{ij}^{pp(k)} l_{ij}^{qq(k)} \hspace{1.9cm}\forall {ij} \in \mathcal{E} 
\end{eqnarray}
\end{small}
\vspace{-0.2cm}

\noindent where, $P_{ij}^{pp(k)},Q_{ij}^{pp(k)},v_{i}^{p(k)},l_{ij}^{pp(k)},l_{ij}^{qq(k)},l_{ij}^{pq(k)}$ are the power flow variables obtained by solving $k^{th}$ iteration of relaxed D-OPF SOCP model. 

The objective is to  gradually reduce the feasibility gap, i.e., $e^{pp(k)}_{ij}$ and $e^{pq(k)}_{ij}$ to sufficiently small values over successive SOCP iterations. This is achieved by enforcing additional directional constraints on the error terms (feasibility gap) defined in (\ref{eq18}), where $\gamma^k < 1$. Notice that, $e^{pp(k)}_{ij} \leq 0$ and $e^{pq(k)}_{ij} \leq 0$. Thus, using (28), the feasibility gap, will be increased towards zero from $k^{th}$ to $(k+1)^{th}$ iteration. Here, $\gamma^{(k)}$ measures the ratio of feasibility gap at $(k+1)^{th}$ and $k^{th}$ iterations  (see \cite{Filter} for details).
\begin{equation}\label{eq18}
    \small
    e^{pp(k+1)}_{ij} \geq \gamma^{(k)} e^{pp(k)}_{ij} \hspace{0.3cm} \text{and} \hspace{0.3cm} e^{pq(k+1)}_{ij} \geq \gamma^{(k)} e^{pq(k)}_{ij}
\end{equation}

Next, upon substituting the expressions for feasibility gaps at $k^{th}$ iteration in (28) and using $\gamma$ where, $\gamma^{(k)} \leq \gamma < 1$, the feasibility gap at $(k+1)^{th}$ SOCP-iteration are expressed as (29) and (30).

\vspace{-0.3cm}
\begin{small}
\begin{eqnarray}\label{eqerror1}
  e^{pp(k+1)}_{ij}  &\geq& \gamma \times \left((P_{ij}^{{pp}(k)})^2 + (Q_{ij}^{pp(k)})^2 - v_i^{p(k)} l_{ij}^{pp(k)}\right)\\
  e^{pq(k+1)}_{ij} &\geq& \gamma \times \left((l_{ij}^{pq(k)})^2  - l_{ij}^{pp(k)} l_{ij}^{qq(k)}\right) 
\end{eqnarray}
\end{small}
\vspace{-0.3cm}

As can be observed  the constraints in (29) and (30) are nonlinear. Thus, we linearize (29) and (30) using first order Taylor series approximation. The resulting linear directional constraints are given in (\ref{direc1})-(\ref{direc2}).
\begin{equation}\label{direc1}
\small
 \begin{split}
 &2P_{ij}^{pp(k)}\Delta P_{ij}^{pp(k+1)}+2Q_{ij}^{pp(k)}\Delta Q_{ij}^{pp(k+1)}-l_{ij}^{pp(k)}\Delta v_i^{p(k+1)} ... \\ 
& ... -v_i^{p(k)}\Delta l_{ij}^{pp(k+1)} \geq  (\gamma-1) e^{pp(k)}_{ij}
\end{split}
\end{equation}
\vspace{-0.3cm}
\begin{equation}\label{direc2}
\small
\begin{split}
&2l_{ij}^{pq(k)}\Delta l_{ij}^{pq(k+1)}-l_{ij}^{pp(k)}\Delta l_{ij}^{qq(k+1)}-l_{ij}^{qq(k)}\Delta l_{ij}^{pp(k+1)} \\ 
&\geq  (\gamma-1) e^{pq(k)}_{ij}
\end{split}
\end{equation}
\noindent where, at $(k+1)^{th}$ iteration, $P_{ij}^{pp(k)}, Q_{ij}^{pp(k)}, v_i^{p(k)}, l_{ij}^{pp(k)}$, $l_{ij}^{qq(k)}, l_{ij}^{pq(k)}, e^{pp(k)}_{ij}, e^{pq(k)}_{ij}$ are known from solving the $k^{th}$ SOCP iteration. Thus, (31) and (32) are linear in unknowns, $\Delta P_{ij}^{pp(k+1)}, \Delta Q_{ij}^{pp(k+1)}, \Delta v_{i}^{p(k+1)},\Delta l_{ij}^{pp(k+1)},$ $\Delta l_{ij}^{qq(k+1)} \text{ and } \Delta l_{ij}^{pq(k+1)}$. 

The actual power flow solution at $(k+1)^{th}$ iteration is obtained by updating power flow variables using (33). 
\begin{equation}
\small
x^{(k+1)} = x^{(k)} + \alpha \Delta x^{(k+1)}
\end{equation}
where, the acceleration factor, $0 < \alpha <1 $, $x^{(k)}$ is the variable obtained at previous iteration, and the change in variables, $\Delta x^{(k+1)}$ is determined at the current iteration.

\subsubsection{Iterative Second-Order Cone Programming (ISOCP)}
The iterative SOCP algorithm for solving three-phase OPF is detailed in this section. We recast the OPF problem in (5)-(16) for simplicity. Let, the problem objective be defined by $f(x)$; the constraints are defined using the following sets of equations: $g_{lin}(x)=0$ and $g_{quad}(x)=0$, where $g_{lin}(x)$ includes (5)-(7) and (11)-(12) (linear in $x$) and $g_{quad}(x)$ includes (8)-(9) (quadratic in $x$). The rest of the constraints, representing operating limits on power flow variables, (13)-(16) are represented as box constraints in (37).

The ISOCP model is obtained by relaxing $g_{quad}(x)=0$ as SOCP constraint, i.e., (24)-(25), and by adding additional directional constraints defined in  (31)-(32), represented as $d(\Delta x)\leq0$ in (38) for simplicity. Thus, $(k+1)^{th}$ iteration of ISOCP based on OPF algorithms is defined as the following. 
\begin{minipage}{23.5em}
		\flushleft
		\small
		 \begin{flalign}\label{eqsocp}
		\text{Minimize:}  \ \ f(x^{k} + \alpha \Delta x^{(k+1)}) &&
		\end{flalign}
		Subject to:
		\vspace{-0.5 cm}
        \begin{flalign}\label{eqsocpc}
        g_{lin}(x^{k} + \alpha \Delta x^{(k+1)})=0 \\
        g_{quad}(x^{k} + \alpha \Delta x^{(k+1)}) \leq 0 \\
        l \leq x^{k} + \alpha \Delta x^{(k+1)} \leq u  \\
        d(\Delta x^{(k+1)}) \leq 0
         \end{flalign}
         \vspace{-0.3cm}
\end{minipage}
where,  $x^{k}$ is known at the $(k+1)^{th}$ iteration. Thus, $(k+1)^{th}$ iteration of ISOCP solves for $\Delta x^{(k+1)}$. 

The proposed algorithm to solve ISOCP model is detailed in Algorithm 2. Same as the PSLP, the problem variables are initialized by solving the three-phase LinDistFlow model \cite{gan2014convex}. At $(k+1)^{th}$ iteration, the ISOCP problem defined in (34)-(38) is solved. Convergence condition is checked using an error variable, ($\varepsilon^k$), defined as the maximum of $|e^{pp(k)}_{ij}|$ and $|e^{pq(k)}_{ij}|$ $\forall {ij} \in \mathcal{E}$. If the $|\varepsilon^k|<tol$, where $tol$ is the pre-specified tolerance, the algorithm stops (\emph{as the current operating point is feasible wrt. original NLP)}; else, the iteration continues. A new operating point at $(k+1)^{th}$ iteration is obtained using (33). 

\textbf{Discussion:} {\em Including a general nonlinear problem objective:} The PSLP algorithm can easily accommodate any nonlinear objective function. Here, we can linearize the objective function along with the constraints (see (19)). However, to formulate the problem as an ISOCP, the objective function needs to be convex. Thus, it would be challenging to work with arbitrary nonlinear objective functions directly. Given that our SOCP model is defined in terms of the square of the voltage magnitude ($v_i^p = (V_i^p)^2$) and branch current flow (${l_{ij}^p} = {{I_{ij}^p}}^2$), the CVR objective (used in this paper), and loss minimization objective can be easily formulated as a linear function in problem variables. However, some approximation will be needed to formulate a general nonlinear objective, such as voltage unbalance factor, as an ISOCP problem. 

\begin{algorithm}[h]
\small
\caption{\small Three-phase OPF using Iterative Second-Order Cone Programming (ISOCP)}\label{alg:ISOCP}
\SetAlgoLined
    \SetKwInOut{Input}{input}
    \SetKwInOut{Output}{output}
\SetKwInput{defi}{define}
\SetKwInput{init}{initialize}
\defi{$x = x_{OPF}$}
\Input{distribution system connectivity model, impedance matrix ($z_{ij}$), load parameters ($p_{L,i}^p$,$q_{L,i}^p$), DG parameters ($p_{DG,i}^p$, $s_{DG,i}^{rated,p}$), operating limit ($I_{ij}^{rated}$), $tol = 0.0001$, $\alpha = 1.0$, $\gamma = 0.9$.}
\Output{Reactive power dispatch from DGs ($q_{DG,i}^p$)}
\init{iteration count, $k = 1$; $x^{(k)} = x_{lin}$, where, $x_{lin}$ is the solution of OPF using LinDistFlow. } \vspace{0.02 cm}
{Calculate: $e_{ij}^{pp(k)}$ and $e_{ij}^{pq(k)}$, $\forall {ij} \in \mathcal{E}$ using (26)-(27)}\\ 
{Calculate: $\varepsilon^{(k)} = \max(|e_{ij}^{pp(k)}|,|e_{ij}^{pq(k)}|)$ $\forall {ij} \in \mathcal{E}$}\\
\While {$|\varepsilon^{(k)}| > tol$}{
{Obtain linearized directional constraints using (31)-(32)}\\
{Solve $(k+1)^{th}$ iteration of SOCP equations in (34)-(38) for $\Delta x^{(k+1)}$} \\
{Update problem variables: $x^{k+1} = x^k + \alpha \Delta x^{(k+1)}$}\\
{Calculate: $e_{ij}^{pp(k+1)},e_{ij}^{pq(k+1)}$, $\forall {ij} \in \mathcal{E}$ using (26)-(27)}\\ 
{Calculate: $\varepsilon^{(k+1)} = \max\left(|e_{ij}^{pp(k+1)}|,|e_{ij}^{pq(k+1)|}\right)$ $\forall {ij} \in \mathcal{E}$} where, $|.|$ is the absolute value operator. \\
{Increment iteration count, $ k = k+1$} \\
}
	\label{algo2}
	\vspace{-0.1 cm}
\end{algorithm}

\subsubsection{Convergence and Optimality of ISOCP Method}
In this section, we detail further analysis of the proposed ISOCP algorithm. Specifically, we prove that the ISOCP iterates converge to the optimal solution of the original nonlinear problem. We define the optimal solution for the original NLP problem as $f^{*}$. The relaxed problem replaces the quadratic equality constraint by SOC constraints. We define $f_{SOCP}^{*}$ as the optimal value of the relaxed problem. However, due to relaxed constraints, $f_{SOCP}^{*}$ is not the actual power flow from the substation. We denote corresponding total power flow from the substation as $f_{Sys}^{*}$.

\noindent \textbf{Lemma 1:} $f_{SOCP}^{*} \leq f^{*}\leq f_{Sys}^{*}$.

\noindent\textbf{Proof:} Since, the feasible domain in the relaxed problem is larger than the original problem, $f_{SOCP}^{*}$ serves as the lower bound for the optimal value of the original problem i.e., $f_{SOCP}^{*} \leq f^{*}$.  Next, let $u_{SOCP}$ is the set of decision variables obtained from solving the relaxed SOCP problem. Then we define $f_{Sys}^{*}$ as the value of the objective function obtained upon implementing the decision variables, $u_{SOCP}$, to the nonlinear power flow model for the given distribution system. Since the SOCP model does not satisfy the original problem constraints, the resulting solution will lead to a different substation power flow, $f_{Sys}^{*}$ from the optimal value of the relaxed problem, $f_{SOCP}^{*}$.  That is, due to the feasibility gap between the relaxed and original nonlinear problem, $f_{socp}^{*}$ and $f_{sys}^{*}$ admit two different value unless the relaxed solution is exact with zero feasibility gap. Here, $f_{Sys}^{*}$ will serve as the upper bound of $f^{*}$. We argue this by contradiction. If $f_{Sys}^{*} \leq f^{*}$, then the decision variables obtained from solving the relaxed problem,  $u_{SOCP}$, leads to the true optimal value for the NLP problem. However, we have previously assumed that the NLP problem's optimal value is $f^{*}$. Therefore, $f^{*} \leq f_{Sys}^{*}$. Thus, $f_{SOCP}^{*} \leq f^{*}\leq f_{Sys}^{*}$. It should be noted that if the discrepancy between the $f_{SOCP}^{*}$ and $f_{Sys}^{*}$ can be reduced, we can achieve the optimal value, $f^{*}$. The iterative procedure proposed in this paper is designed specifically to achieve this goal.

Next, we define error functions:  $e_{ij}^{pp} = (P_{ij}^{pp})^2 + (Q_{ij}^{pp})^2 - v_i^p  l_{ij}^{pp}$, and $e_{ij}^{pq} = (l_{ij}^{pq})^2 - l_{ij}^{pp}  l_{ij}^{qq}$. Let, $e(x)=\left[\begin{array}{c}
e_{ij}^{pp}  \\
e_{ij}^{pq} \end{array} \right] 
$. Then, the proposed iterative approach aims at gradually reducing $|e(x^{(k)})|$ to zero as $k\rightarrow \infty$, where $e(x^{(k)})$ are set of errors at $k^{th}$ iteration of the ISOCP.

\noindent\textbf{Lemma 2:} $f^{*}$ can be achieved if $\lim_{k\rightarrow \infty} |e(x^{(k)})| = 0$.

\textbf{Proof:} Let, $f_{SOCP}^{(k)}$ be the optimal value of the relaxed (SOCP) problem solved at $k^{th}$ iteration and $f_{Sys}^{(k)}$ be the corresponding substation power flow. Then, if the SOCP solutions satisfy the original power flow constraints, then, $e(x^{(k)})=0$ and $f_{SOCP} = f_{Sys}$. In other words, if $\lim_{k\rightarrow \infty} |e(x^{(k)})| = 0$, $\lim_{k\rightarrow \infty} \left(f_{Sys}^{(k)}-f_{SOCP}^{(k)}\right) = 0$. Further, using Lemma 1, $f_{SOCP}^{(k)} \leq f^{*}\leq f_{Sys}^{(k)}$, if $\lim_{k\rightarrow \infty} |e(x^{(k)})| = 0$, the actual optimal value for the problem, i.e., $f^{*}$ is achieved as $k \rightarrow \infty$. Here, we construct iterates to specifically satisfy the property, $\lim_{k\rightarrow \infty} |e(x^{(k)})| = 0$.

\textbf{Theorem 1:} Let, $\gamma^{(k)} = \dfrac{|e(x^{(k)})|}{|e(x^{(k-1)})|}$. Then, at $k^{th}$ iteration, if we impose constraint, $\gamma^{(k)}\leq \gamma < 1$ to the relaxed (SOCP) problem, then the resulting iterative problem will converge to the optimal solution, $f^{*}$, as $k\rightarrow \infty$.

\noindent\textbf{Proof:} Imposing constraint $\gamma^{(k)}\leq \gamma < 1$ where, $\gamma^{(k)} = \dfrac{|e(x^{(k)})|}{|e(x^{(k-1)})|}$ is equivalent to:
\begin{equation}
e(x^{(k)}) \geq \gamma e(x^{(k-1)})
\end{equation}
On applying the above inequality recursively, we obtain: 
\begin{equation}
e(x^{(k)}) \geq \gamma^k e(x^{(0)})
\end{equation}
Since, $\gamma < 1$, as $k\rightarrow \infty$, $\gamma^k e(x^{(0)}) \rightarrow 0$. Since, $ e(x^{(k)}) \leq 0$, we can conclude,  $\lim_{k\rightarrow \infty} |e(x^{(k)})| = 0$. Thus, upon incorporating the constraints of the form (39) in relaxed (SOCP) problem, we can ensure that $\lim_{k\rightarrow \infty} |e(x^{(k)})| = 0$ and consequently $\lim_{k\rightarrow \infty} \left(f_{Sys}^{(k)}-f_{SOCP}^{(k)}\right) = 0$. Together with the fact that $f_{SOCP}^{*} \leq f^{*}\leq f_{Sys}^{*}$, the optimal solution to the original problem, $f^{*}$, is achieved as $k\rightarrow \infty$.

\textbf{Linearizing the inequality constraints and further convergence arguments:}
Unfortunately, the inequality constraint (39) is non-convex. Thus, we use a first order taylor series approximation to linearize $e(x^{(k)})$ around the solutions from the previous iterate, i.e. $x^{(k-1)}$. This results in linear equations of the form (41). 
\begin{equation}
{\triangledown_{x^{(k-1)}}e(x^{(k-1)})}\Delta x^{(k)} \geq (\gamma-1)e(x^{(k-1))}
\end{equation}

Notice that since, $e(x^{(k)})$ is a convex quadratic function, the first-order Taylor series approximation will be the global under-estimator of $e(x^{(k)})$. Therefore, upon imposing linearized constraints (39), inequality $e(x^{(k)}) \geq \gamma e(x^{(k-1)})$ can be ensured. Thus, the convergence of the iterative procedure can always be ensured. However, since the search direction is restricted to half-space in (41) as opposed to non-convex set in (39), we cannot guarantee global optimality upon including constraint (41) instead of (39) in the SOCP formulation. 

 It should be noted that at given iteration $(k)$, $x^{(k-1)}$ and hence the value of $e(x^{(k-1)})$ are known. Thus, (39) is a set of inequality constraints in the problem variables at the current iteration $x^{(k)}$. That is, $e(x^{(k)}) \geq C$, where, $C = \gamma e(x^{(k-1)})$ is a constant vector where, each element of $C$ is $\leq0$, by definition. At the optimal solution for the original NLP, by construction, $e(x_{opt}) = 0$ (i.e. feasibility gap is zero). Then, at any given iteration, $e(x^{(k)}) \geq C$ will always admit a solution. Thus, (39) and its linear approximation (41) will admit solutions at a given iteration.

\section{Results and Discussions}
The proposed three-phase D-OPF formulations are validated using two distribution test feeders: (1) IEEE 123-bus\cite{IEEE123} (with 267 single-phase nodes) and (2) modified R3-12.47-2 feeder (with 860 single-phase nodes). Multiple test cases are simulated to evaluate the feasibility, optimality, and computation time for the two algorithms against the equivalent NLP D-OPF model. The NLP D-OPF is solved using a commercial solver, `Knitro' interfaced with MATLAB. The proposed PSLP and ISOCP algorithms are implemented using MATLAB and solved using CPLEX 12.7 and Gurobi, respectively.

\vspace{-0.4cm}     
\subsection{Verification of Approximate BFM Formulations in (5)-(9)}
In this section, the power flow solutions from approximate BFM model in (5)-(9) are validated against those obtained using OpenDSS. Since the proposed D-OPF algorithms employ approximate BFM, this validation step is crucial. The largest errors in apparent power flow and bus voltages are reported for the two test feeders in Table I. 

It is observed that for the IEEE-123 bus system, the largest errors in power flow and node voltages are 0.162\% and 0.0025 (pu). Similarly, for the modified R3-12.47-2 feeder, the largest errors in power flow and node voltages are observed to be 0.068\% and 0.0007(pu), respectively. It should be noted that the errors in power flow variables are very small. Thus, the proposed approximate BFM accurately represents an unbalanced distribution system. 

\begin{table}[h]
		\centering
		\vspace{-0.4cm}       
		\caption{Approximate NLP BFM Model in (5)-(9) vs. OpenDSS}
     \vspace{-0.2cm}      
		\label{singletable}
		\begin{tabular}{c|c|c|c}
		    \toprule[0.4 mm]
			\hline
			{Test Feeder} & \% Loading & $S_{flow}(\%)$ & $V$(pu.)\\
			\hline
            \hline
            {IEEE 123 Bus}& 100\%& 0.162 & 0.0025 \\
            \hline
            {R3-12.47-2}& 100\%& 0.068 & 0.0007\\
			\toprule[0.4 mm]
		\end{tabular}
\end{table}

Note that the approximate BFM is obtained by assuming (1) node voltage phase angles are $120^0$ apart and (2) branch current phase angles can be approximated using the solution for an equivalent constant impedance load model. We further validate these assumptions using the two test feeders. It is observed that the largest deviation in the node voltage phase angle difference from $120^0$ is around $2^0$. Also, the largest deviation in the phase current angle from that obtained using a constant impedance load model is around $2^0$. Thus, the assumptions are reasonably accurate for unbalanced distribution systems. Finally, it should be noted that we have presented an extensive validation of the proposed approximate BFM model in our prior work under several different loading conditions with varying levels of voltage unbalance \cite{jha2019bi}. However, we have included a few scenarios here for completeness.

\vspace{-0.2 cm}

\subsection{Benchmark NLP Algorithms}
In this section, we benchmark different NLP algorithms to solve the NLP model for the D-OPF problem using `Knitro'. The solver includes three different NLP algorithms: (1) Interior-point algorithm, (2) Sequential Quadratic Programming (SQP) algorithm, (3) Active Set algorithm (where it uses a sequential linear-quadratic programming (SLQP) algorithm, similar in nature to a sequential quadratic programming method but using linear programming sub problems to estimate the active set) \cite{andrei2017active}. In what follows, we detail the convergence results for the 123-bus test system using the above three NLP algorithms. The results are summarized in the Table II. We have also compared the resulting optimality and feasibility errors used as the termination criteria in `Knitro’ solver. The optimality error is defined as the maximum violation of the first-order conditions for identifying a locally optimal solution at the current iterate. Feasibility error is defined as the maximum violation of problem constraints at the current iterate.  As can be seen, interior-point method resulted in high values of feasibility and optimality errors and failed to converge. While both SQP and Active-set methods were able to reduce the feasibility errors for all test cases, SQP led to a higher value of  optimality errors and converged to a higher value of the objective function. Thus, owing to better convergence properties, we adopt Active-set algorithm in `Knitro' to solve the NLP problem for both 123-bus and modified R3-12.47-2 test feeders.

The NLP model requires a good warm start for better performance. Here, we warm start the NLP OPF using the linearized AC OPF model's solutions for the same objective function \cite{gan2014convex}. Note that the linear AC power flow model very well approximates the node voltages and provides good initialization for the NLP model for distribution OPF problems.

\begin{table}[t]
\centering
\small
\caption{Comparison of different NLP D-OPF algorithms for IEEE 123 bus system}
\label{NLPsolvecom}
\begin{tabular}{c|c|c|c|c}
\hline
\hline
\% DG & Without D-OPF & Interior Point & SQP & Active set \\
\hline
&\multicolumn{4}{c}{Substation Power (MW)}\\
\hline
 10 & 3.329 & infeasible &  3.265 & 3.264 \\  
 \hline
  20&  2.913 & infeasible & 2.863 & 2.848  \\
 \hline
  30& 2.502 & infeasible & 2.458& 2.431 \\
  \hline
  40&2.098 & infeasible & 2.05 &2.02\\
  \hline
  50&1.903 & infeasible & 1.822&1.820 \\
  \hline
   &\multicolumn{4}{c}{Optimality Error}\\
   \hline
 10 & - & $8.40e^3$ & $5.04e^{-4}$ & $4.62e^{-11}$ \\
 \hline
  20 & -& 53.8 &$2.59e^{-2}$& $2.61^{-12}$ \\
 \hline
  30& - &1.01 & $1.12e^{-2}$& $6.36e^{-13}$\\
  \hline
  40&- & 1.10 & $1.10e^{-2}$ & $2.61^{-11}$\\
  \hline
  50&- &1.06 & $2.09e^{-2}$ & $7.85e^{-14}$\\
  \hline
   &\multicolumn{4}{c}{Feasibility Error}\\
   \hline
 10 & - & $8.41e^3$ & $8.42e^{-9}$ & $2.63e^{-11}$\\  
 \hline
  20&- & 0.991 & $1.31e^{-9}$ & $2.35e^{-9}$ \\
 \hline
  30& -& 1.31& $2.51e^{-9}$ & $4.54e^{-9}$\\
  \hline
  40&- & 1.73 & $4.32e^{-9}$ &$ 2.06^{-10}$\\
  \hline
  50&- &10.0 & $6.49^{-9}$ & $2.65e^{-11}$ \\
 \hline
\end{tabular}
\vspace{-0.2 cm}
\end{table}

\textbf{Discussion - }{\em Convergence of NLP Problem:} The `Knitro' solver uses derivative check to ensure that the solution has reached a local optimum value. The optimality for the NLP solutions is ensured using two error parameters in the `Knitro' solver: (1) optimality error, and (2) feasibility error. Specifically, the nonlinear optimizer in the `Knitro' uses optimality error parameter to check the convergence to the local optimum solution and feasibility error parameter to ensure that constraints are satisfied at the current iterate. For the NLP model, we have used the following thresholds for convergence in `Knitro’: optimality error tolerance  $= 10^{-6}$, feasibility error tolerance  $= 10^{-6}$.

\subsection{Performance of the Proposed D-OPF Algorithms}
The proposed PSLP (approximation-based) and ISOCP (relaxation-based) D-OPF algorithms are evaluated for the optimality, feasibility and compute time using the two test feeders. The D-OPF solutions obtained using the proposed algorithms are compared against those obtained upon solving the actual NLP D-OPF problem. We also compared the compute time required to solve the proposed D-OPF algorithms with respect to the NLP D-OPF model. Finally, the D-OPF solutions are verified using a distribution system simulator, OpenDSS. For all simulations, we have assumed loads to have a CVR factor of 2.0. Note that the CVR values are arbitrary and can be easily adjusted based on the parameters for the ZIP model of the load, if available (see \cite{jha2019bi}).

\begin{table}[t]
		\centering
		\caption{Comparison of D-OPF algorithms for IEEE 123 bus system}
		\vspace{-0.3cm}
		\label{IEEE123buscomp}
		\begin{tabular}{L|c|c|c|c|L}
		\hline
		\hline
		\% DG & Without D-OPF & NLP & PSLP &  ISOCP &  \multirow{2}{=}{OpenDSS Validation} \\
        \cline{2-5}
        &\multicolumn{4}{c|}{Substation Power (MW)}&\\
        \hline
         10& 3.329 & 3.264 & 3.268 & 3.261 & 3.262 \\
 		\hline
 		 20&2.913 & 2.845 & 2.848 & 2.842 & 2.843 \\
 		\hline
 		 30& 2.502 & 2.431 & 2.433 & 2.426 & 2.428 \\
 		 \hline
 		 40& 2.098 & 2.022 & 2.024 & 2.018 & 2.018 \\
 		 \hline
 		 50& 1.903 & 1.820 & 1.821 & 1.814 & 1.816 \\
 		\hline
       &&\multicolumn{3}{c|}{Computation time (secs)}&\\
        \hline
         10& - & 342 & 3.9 & 98  &  -\\
        \hline
 		 20&- & 340 & 4.2 & 91 & - \\
 		\hline
 		 30& -& 338 & 7.1 & 77 & - \\
 		 \hline
 		 40& -& 299 & 8.3 & 70 & - \\
 		 \hline
 		50& -& 219 & 12.2 & 63 &-  \\
 		\hline
 		\hline
\end{tabular}
\end{table}

\subsubsection{IEEE 123-bus test system}
The IEEE 123-bus system is a medium-sized feeder (with 267 single-phase nodes) with unbalanced loading and several single-phase lines and loads. The feeder is modified to include DGs with smart inverters. Multiple test cases are simulated to show the effect of increasing the number of control variables (number of DGs with smart inverters) on the performance of the proposed D-OPF algorithms. Specifically, we vary the percentage of DG penetrations from 10\% to 50\% in the steps of 10\%.  The increase in DG penetrations implies that new controllable DG units are added to the distribution system. Each smart inverter is rated at 48 kVA. 

Recall, the problem objective is CVR that is to minimize the active power consumption from the substation by coordinating the reactive power dispatch from the smart inverters. Here, the D-OPF problem detailed in (5)-(16) is solved for different DG penetration levels. The results obtained upon solving different D-OPF algorithms, viz. NLP, PSLP and ISOCP is shown in Table. \ref{IEEE123buscomp}.  It provides the substation power demand with and without D-OPF for different levels of DG penetration. It can be observed that as the DG penetration is increased, the percentage reduction in substation power, upon solving D-OPF for CVR objective, is also increased. This is because, at higher DG penetrations, additional reactive power support is available from the smart inveretrs that can be utilized to regulate feeder voltages towards their lower allowable limit in order to extract additional CVR benefits. 

\begin{figure}[t]
 \centering
  \vspace{-0.3cm}
\includegraphics[width=3.2in]{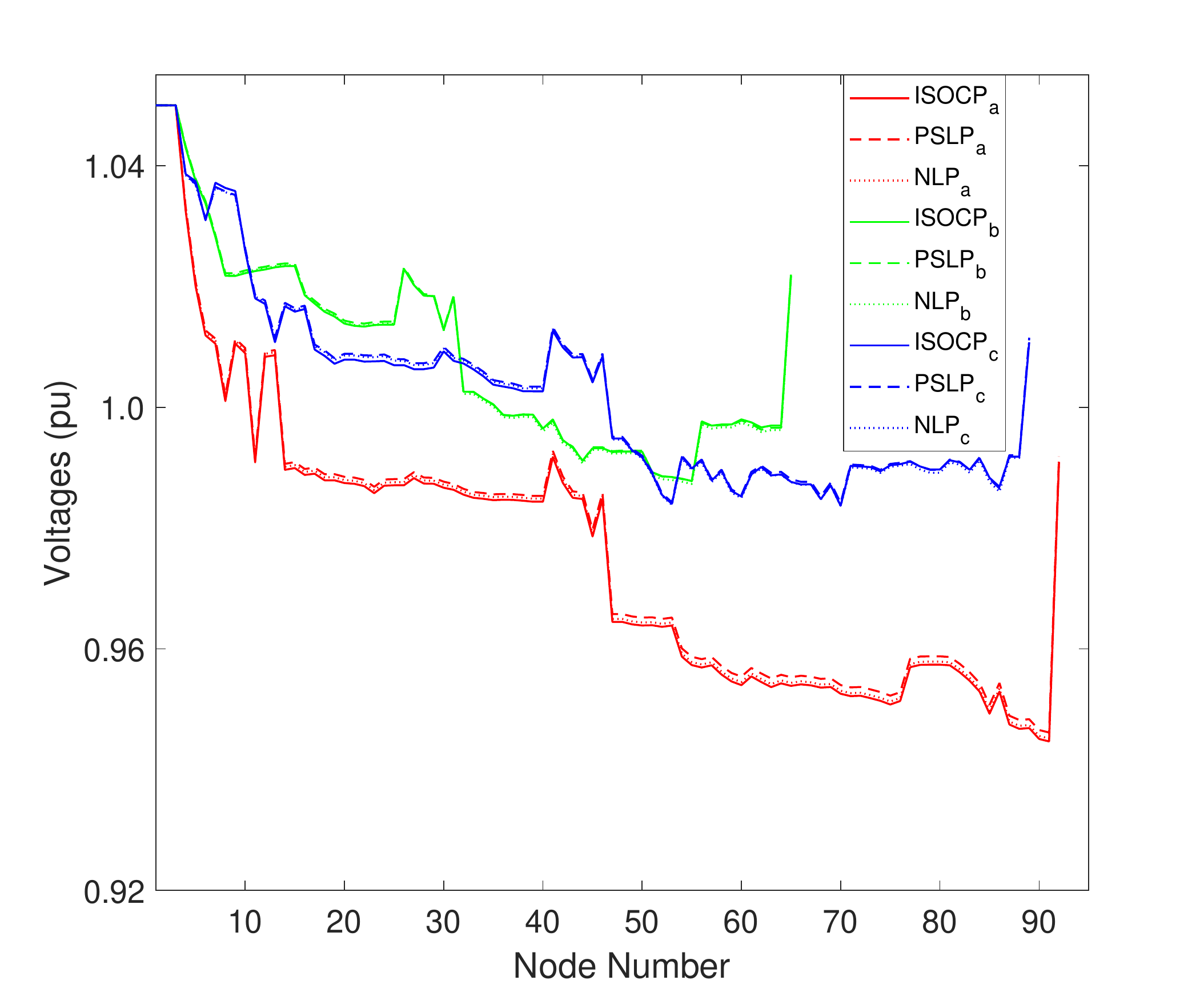}
	\caption{Phase voltage with 50\% DG penetration for IEEE-123 node system}
	\label{phasevoltage}
\end{figure}

Next, we present the comparative analysis of the proposed D-OPF algorithms. Table. \ref{IEEE123buscomp} shows the optimal substation power demand upon solving the proposed PSLP and ISOCP D-OPF algorithms that are based on the approximation and relaxation of the NLP D-OPF model, respectively. It should be noted that the optimum values from the proposed D-OPF algorithms closely matches to those obtained by solving the original NLP D-OPF model. Further, the simulation results from OpenDSS upon implementing the decision variables obtained from D-OPF models also closely matches with the D-OPF solutions.  For results in Table III, we implement the decision variables obtained from solving the NLP model to OpenDSS. Further, we also implement the decision variables obtained from all three OPF models to OpenDSS and observe bus voltages for three test cases.   The voltages profile for the 50\% DG penetration case for IEEE 123-bus test system is shown in Fig.\ref{phasevoltage}. Note that the node voltages vary all the way from 1.05 pu to 0.95 pu (the acceptable range for the distribution voltages), thus showing the presence of important nonlinearities in the problem formulation. Further, it can be observed that the voltage profile obtained from the three D-OPF algorithms closely match, thus validating that the solutions of the proposed algorithms converge to the same solution point. In short, the simulation results validate the following: (1) the proposed PSLP and ISOCP D-OPF algorithm arrive at the same optimal results as those obtained via original NLP D-OPF; and (2) the power flow solutions obtained from the three D-OPF models are feasible for the given unbalanced distribution system.

Recall that the SOCP relaxation was inexact and we proposed an iterative approach (with a directional constraint) to reach towards a feasible power flow solution over successive iterations. To further elaborate the feasibility aspects of ISOCP algorithm, we plot the feasibility gap, $e^{k} =  \max\left(|e_{ij}^{pp(k+1)}|,|e_{ij}^{pq(k+1)|}\right)$, for each DG penetration level (see Fig. 2).  We have used the power base of 1MW for per-unit conversion. Based on the 1MW base, a feasibility gap of 1.5 pu corresponds to 1.5MW$^2$ of feasibility gap at the substation bus.  It can be observed that the feasibility gap decreases over multiple iterations of SOCP and converges to a feasible OPF solution. Finally, the computational advantages of the proposed algorithms are apparent from the Table III. While it takes close to 6 min to solve the NLP D-OPF, the PSLP and ISOCP solves within 10 sec and 1.5 mins, respectively. Note that each iteration of SOCP take approximately 7 sec to solve.

\begin{figure}[t]
\includegraphics[width=3.6in]{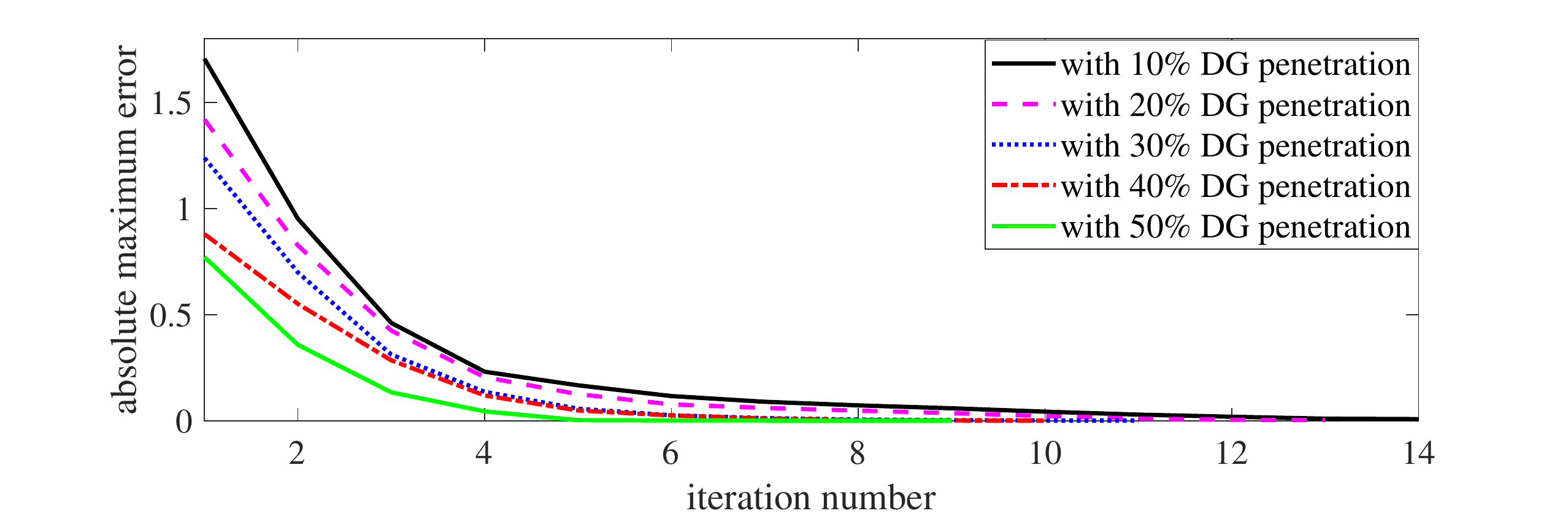}
\vspace{-0.5cm}
	\caption{IEEE-123 bus: Reduction in feasibility gap vs. number of iterations.}
	\label{feasvoil123bus}
\end{figure}

\begin{table*}[t]
		\centering
		\caption{Comparison of D-OPF algorithms for IEEE 123 bus system }
		\label{IEEE123buscomp1}
		\begin{tabular}{c|c|c|c|c|c}
		\hline
		\hline
		CVR factor and $S_{DG}$ rating & \multicolumn{4}{c|}{Substation Power (MW)}& \multirow{2}{*}{CVR benefits (kW)} \\ 
		\cline{2-5}
		&Without D-OPF  & NLP & PSLP &  ISOCP &  \\
 		\hline
 		CVR = 2  and $S_{DG} = 1.2 P_{rated}$ &2.885& 2.798& 2.798 & 2.794 & 87 \\
 		\hline
 		CVR = 3  and $S_{DG} = 1.2 P_{rated}$ & 2.905& 2.746 & 2.748 &2.741 & 159 \\
 		 \hline
 		CVR = 2  and $S_{DG} = 1.3 P_{rated}$ &2.885& 2.794& 2.796 & 2.790 & 91 \\
 		\hline
 		CVR = 3  and $S_{DG}= 1.3 P_{rated}$ & 2.905& 2.735 & 2.738 & 2.732& 170 \\
 		 \hline
 		\hline
 		 \end{tabular}
\end{table*}

We have provided additional simulation cases to further validate the accuracy of proposed scalable algorithms. Here, we show the effects of the available reactive support and CVR factors for the load, on CVR objective (see Table \ref{IEEE123buscomp1}). The simulations correspond to 40\% DG penetration case (at the rated loading condition) but with PV generating active power at 50\% of its rated capacity. This increases the available reactive power support from the inverters. Further, we create cases with rated PV capacity of 120\% and 130\% of the active power rating. We also vary the load CVR factor while keeping the reactive power support constant. The results are shown in Table \ref{IEEE123buscomp1}. As can be observed, the CVR benefits (measured as the reduction in substation power demand) increases upon increasing the available reactive power support and the load CVR factor. Further, all three algorithms converge to same solution for all test cases. Notice that the similarity of the optimal solutions is argued from an engineering perspective.

\begin{figure}[t]
 \centering
  \vspace{-0.3cm}
\includegraphics[width=3.2in]{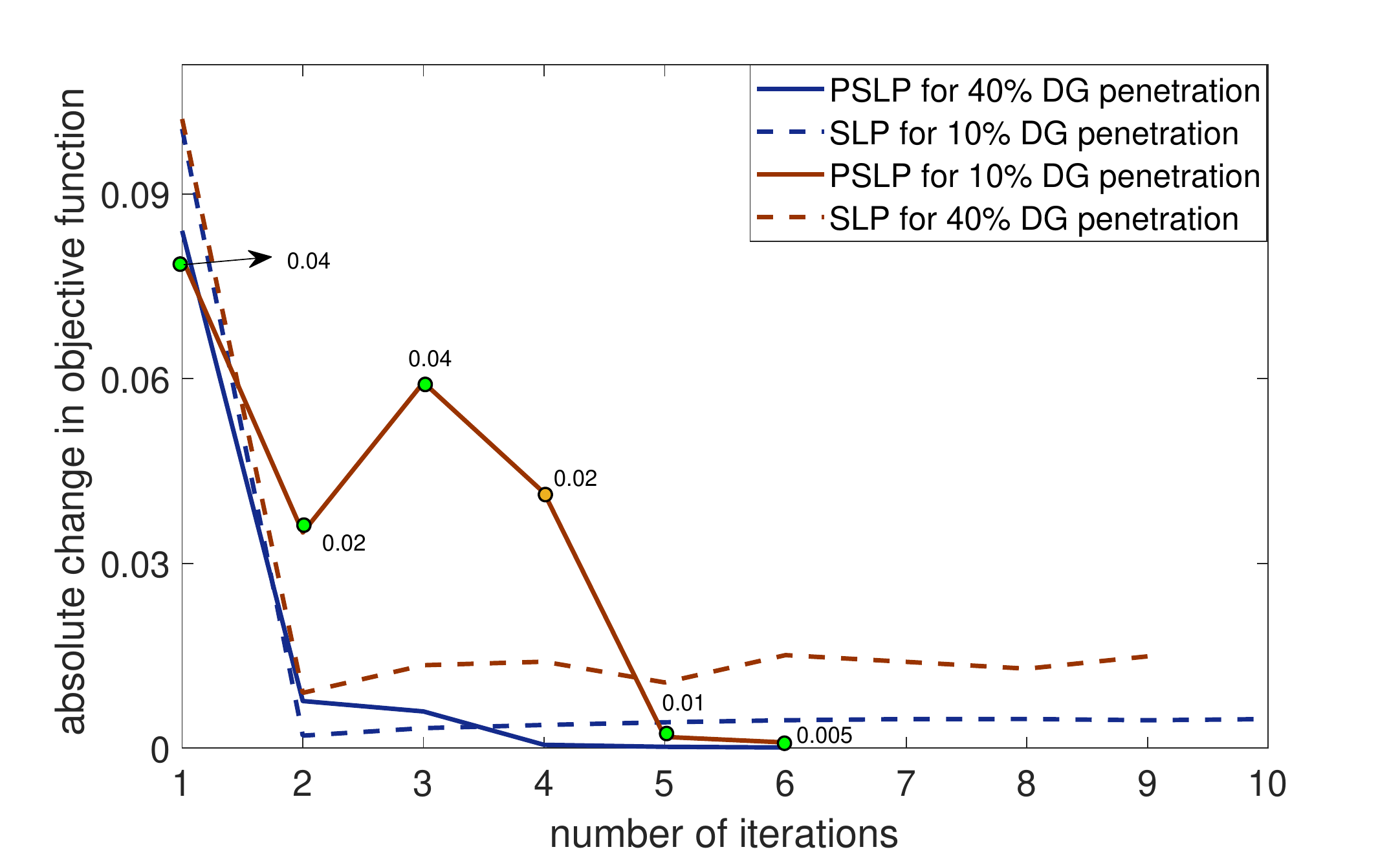}
 \vspace{-0.3cm}
	\caption{The performance comparison of the SLP and PSLP algorithm}
	\label{SLPobjective}
\end{figure}

It is to be noted that successive linear programming (SLP) is a popular linearization approach to solve the NLP problems. In SLP, the nonlinear problem is linearized around the current operating point at every iteration, and the resulting LP is solved to update the problem variables. The process repeats until a converged solution is obtained. The convergence is typically identified when the change in consecutive iterates is within the pre-specified tolerance. While SLP is efficient and fast in solving nonlinear problems, it is known to result in convergence issues  \cite{TEdgar}. The SLP algorithm's primary challenge is that the solutions oscillate as they get closer to the optimal value and may not converge or take several additional steps to converge to the optimal solution. These oscillations are due to the fixed value of the step bounds \cite{TEdgar}. In related literature from the optimization community, the SLP algorithms have been improved by using a penalty function where adaptive step bounds are enforced on the control variables. The penalty function and adaptive step bounds lead to more superior convergence properties \cite{TEdgar}. For comparison, we simulated the SLP and the PSLP algorithms for the IEEE 123-bus test system for 10\% and 40\% DG penetration. We show the plot for the absolute value of the error term for the consecutive iterates in  Fig.~\ref{SLPobjective}. It can be observed from the figure that for the SLP, the error term oscillates and does not show proper convergence. Contrary to SLP, in the PSLP, the step bounds are adaptive; that is, they are reduced as the solution progressively converges to the optimal solution. Notice that for 10\% penetration case, while the error increases at iterate 3, with proper setting of the step bounds, the PSLP algorithm converges within the next two iterations. However, the SLP algorithm could not converge for the simulated test cases.     

\subsubsection{R3-12.47-2 test system} 
Next, we demonstrate the proposed D-OPF algorithms using a larger test feeder, R3-12.47-2. The selected test feeder includes 329 buses with a total of 860 single-phase nodes. Here, the total number of power flow variables are 4225. The DG penetration levels are varied from 10\%-50\% in the steps of 10\% to evaluate the effect of the number of control variables on the performance of the D-OPF algorithms. Each smart inverter is rated at 240 kVA. The results are shown in Table \ref{329buscomp}.

As observed for the IEEE 123-bus test system, both PSLP and ISOCP algorithms lead to same optimal solution as the NLP D-OPF problem for this test system. Further, the power flow results obtained from the D-OPF models are feasible for the unbalance power distribution system as validated using OpenDSS.  For  results  in  Table  V,  we  implement  the decision  variables  obtained  from  solving  the  NLP  model  to OpenDSS. The ISOCP algorithm is successfully able to reduce the feasibility gap over successive SOCP iterations (see Fig. 4). Further, the D-OPF is successfully able to reduce the feeder power consumption for all DG penetration levels where the percentage reduction in demand increases with the DG penetrations. 

For the larger test feeder, the greatest advantage of using proposed approximate and relaxed D-OPF models is evident when comparing the compute time. The NLP D-OPF algorithm takes close to 45-60 mins to obtain a converged solution for the selected test feeder. While the PSLP algorithm takes less than 1 min to obtain a converged solution. As the PSLP based algorithm solves a linear sets of equations, the computation time required to solve the problem is significantly less. Similarly, the  ISOCP algorithm is also significantly faster as the maximum time   to reach a converged solution is 6 mins. In this case, every iteration of ISOCP solution takes 20 sec and the feasibility gap is reduced below the tolerance in around 10-20 iterations. It should be noted that the number of iteration required to reduce the feasibility gap (absolute maximum error) depends on the initial gap.  While ISOCP takes longer to solve, it (1) includes the entire feasible space of the original NLP problem thus theoretically can reach same optimization solution as the NLP problem; (2) is more stable than the PSLP approach that is prone to oscillations and poses convergence issues when the linear approximations are not valid for the selected trust region. 

\noindent \textbf{Discussion - }{\em Voltage Control using Legacy and New Devices:} Notice that a small saving is observed for the OPF problem in Tables III, IV and V. This is because only the reactive power support from the DGs is utilized for voltage control. Typically, the kVA rating of the smart inverters is 120\% or 130\% of their peak active power generation capacity. Thus, they have relatively small reactive power support available for voltage control. Typically, in a distribution system, both legacy devices (capacitor banks and voltage regulators) and the smart inverters are utilized for voltage control to realize the CVR benefits. The combined management of these devices leads to higher CVR benefits, as observed in our previous paper \cite{jha2019bi}.   

\begin{table}[t]
		\centering
		\caption{Comparison of D-OPF algorithms for R3-12.47-2 system}
		\vspace{-0.3cm}
		\label{329buscomp}
		\begin{tabular}{c|L|c|c|c|L}
		\hline
		\hline
		\% DG & Without D-OPF & NLP & PSLP &  ISOCP &  \multirow{2}{=}{OpenDSS Validation} \\
        \cline{2-5}
        &\multicolumn{4}{c|}{Substation Power (MW)}&\\
        \hline
         10 & 3.775 & 3.749 & 3.750 & 3.748 & 3.749 \\
         \hline
         20 & 3.328 & 3.278 & 3.281 & 3.275 & 3.277 \\
 		\hline
 		 30& 3.031 & 2.965 & 2.966 & 2.963 & 2.964 \\
 		 \hline
 		 40& 2.439 & 2.342 & 2.345 & 2.341 & 2.341 \\
 		 \hline
 		 50& 1.847 & 1.737 & 1.739 & 1.736 & 1.737 \\
 		\hline
        &  &\multicolumn{3}{c|}{Computation time (secs)}&\\
        \hline
         10& - & 2531 & 22.2 & 378 &-  \\
         \hline
         20& - & 2865 & 32.2 & 306 &- \\
 		\hline
 		 30&- & 2930 & 47.8 & 234 &- \\
 		 \hline
         40& - & 3413 & 36.2 & 198 &- \\
 		 \hline
 		 50&- & 3381 & 41.2 & 180 &- \\
 		\hline
 		\hline
\end{tabular}
\vspace{-0.2 cm}
\end{table}

\begin{figure}[t]
\includegraphics[width=3.6in]{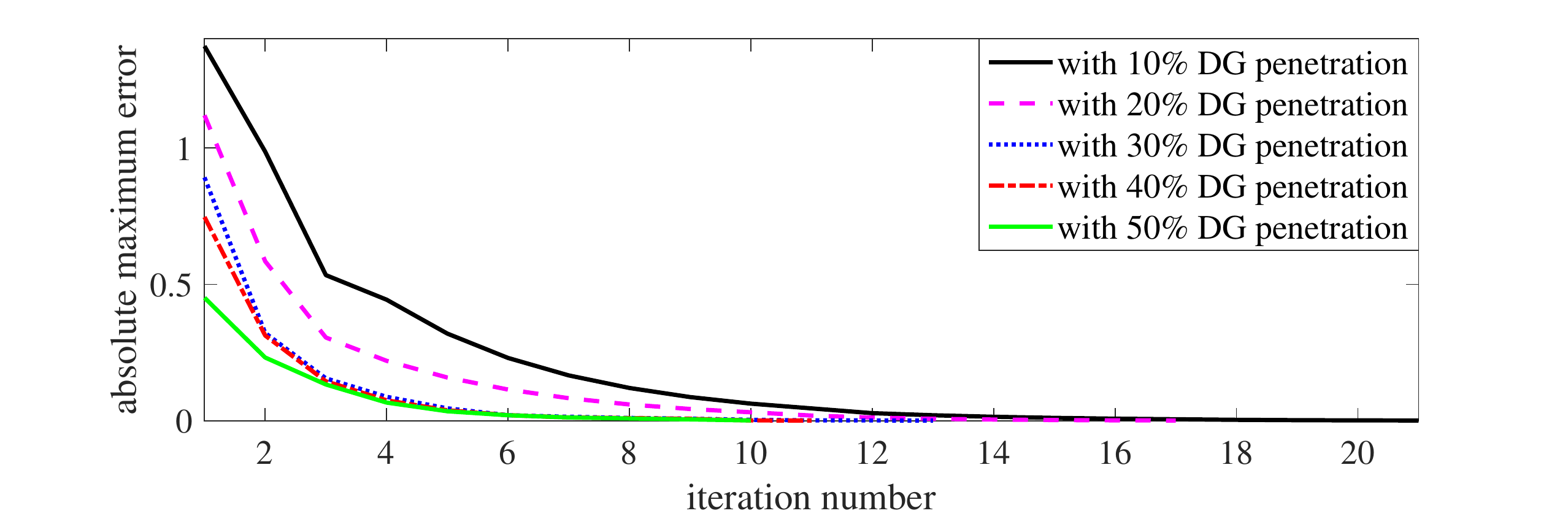}
\vspace{-0.8cm}
	\caption{R3-12.47-2 test system: Reduction in feasibility gap with iterations.}
\vspace{-0.3 cm}
	\label{feasvoil329bus}
\end{figure}

\section{Conclusion}\label{Conclusion}
The distribution-level optimal power flow (D-OPF) problems are specifically challenging due to added nonlinearities posing convergence and scalability issues when applied to a large unbalanced distribution system. The existing scalable D-OPF algorithms either lead to suboptimal solutions due to invalid approximation or infeasible solutions upon solving relaxed problems that are inexact for the unbalanced power distribution systems. In this paper, we present two novel D-OPF algorithms based on approximation and relaxation techniques for unbalanced power distribution systems that specifically address the aforementioned gap in the D-OPF literature. The approximation leads to a PSLP problem that solves successive linear approximations of the NLP D-OPF problem. The relaxation leads to an ISOCP problem that iteratively solves SOCP problems with directional constraints to obtain an optimal and feasible D-OPF solution. Here, our motivation is to develop iterative algorithms that solve simpler D-OPF sub-problems using approximate/relaxed models, but are simultaneously able to reach to the solution of the original NLP over successive iterations. Each iteration of the proposed D-OPF subproblem is of much lower complexity than the original NLP D-OPF. Further, both D-OPF algorithms successfully converge to the solution of the original NLP D-OPF problem while significantly decreasing the compute time. 

Notice that both methods, approximation and relaxation, have their pros and cons. While the LP problem obtained by approximation is much faster than the SOCP (obtained by relaxation), it requires tuning additional parameters wrt. the trust-region of update variables. On the other hand, since, the relaxed model includes the entire feasible space of the original NLP problem, it should be able to reach the same solution as the NLP problem. Further work is needed to understand the tradeoffs between the two techniques for the D-OPF problem.

\bibliographystyle{IEEEtran}
\bibliography{references}

\end{document}